%Paper: funct-an/9410004
%From: Roland Speicher <L95@VM.URZ.UNI-HEIDELBERG.DE>
%Date: Mon, 17 Oct 94 13:18:24 CET

\input amstex \documentstyle {amsppt}
\magnification=1200
\hsize=15truecm
\baselineskip=12pt
\hoffset=0truecm
\predefine\lll{\l}
\redefine\cA{{\Cal A}}
\redefine\cM{{\Cal M}}

\redefine\cV{\pi}

\redefine\RR{\text{{\bf R}}}
\redefine\NN{\text{{\bf N}}}
\redefine\ZZ{\text{{\bf Z}}}
\redefine\CC{\text{{\bf C}}}

\redefine\ff{\varphi}

\redefine\cA{{\Cal A}}

\redefine\lb{\lbrack}
\redefine\rb{\rbrack}
\redefine\l{$\lbrack$}
\redefine\r{$\rbrack$}

\redefine\wlim{\operatornamewithlimits{w-lim}}

\redefine\cA{{\Cal A}}

\redefine\sqp{\sqcup \hskip -0.9em \sqcap \hskip -0.95em +}

\redefine\la{\langle}
\redefine\ra{\rangle}

\redefine\splus{\sqp}
\redefine\ker{\text{ker}}
\redefine\cupdot{\cup}
\redefine\hateq{\hat =}
\redefine\NCzn{NC_2(2n)}
\redefine\mh{\hat\mu}
\redefine\mN{\mu_N}
\redefine\nN{\nu_N}
\redefine\nh{\hat\nu}
\redefine\amb{(\alpha^2-\beta^2)}

\redefine\ak{a^{(k)}}
\redefine\yk{y^{(k)}}
\redefine\yks{y^{(k)*}}
\redefine\aks{a^{(k)*}}
\redefine\al{a^{(l)}}
\redefine\yl{y^{(l)}}
\document
\heading
{\bf CONVOLUTION AND LIMIT THEOREMS FOR}\\
{\bf CONDITIONALLY FREE RANDOM VARIABLES} \endheading
\heading
{\bf Marek Bo\D zejko$^1$, Michael Leinert$^2$, Roland Speicher$^2$}\\
\quad\\
$^1$ Instytut Matematyczny\\
Uniwersytet Wroc\lll awski\\ Plac Grunwaldzki 2/4\\
50-384 Wroc\lll aw\\ Poland \\ e-mail: bozejko\@math.uni.wroc.pl\\
\quad\\ and\\ \quad\\
$^2$ Institut f\"ur Angewandte Mathematik\\ Universit\"at Heidelberg\\
Im Neuenheimer Feld 294\\ D-69120 Heidelberg\\ Germany\\
e-mail: C21\@vm.urz.uni-heidelberg.de (M.L.)\\
\quad\qquad L95\@vm.urz.uni-heidelberg.de (R.S.)\\
\quad\\
{\bf Dedicated to Professor Wilhelm von Waldenfels}\\
{\bf on the occasion of his 60th birthday}
\endheading
\vskip1cm
\heading
{\bf Abstract}
\endheading
We introduce the notion of a
conditionally free product and conditionally free convolution.
We describe this convolution both from a combinatorial point of view,
by showing its connection with the lattice of non-crossing partitions, and
from an analytic point of view, by presenting the basic formula for its
$R$-transform. We calculate
explicitly the distributions of the conditionally free
Gaussian and conditionally free Poisson distribution.
\newline\par\vskip1cm
\pagebreak
\heading
{\bf 1. Introduction}
\endheading
In \l BSp\r, we introduced a generalization
with respect to two states
of the reduced free product
of Voiculescu \l Voi1,VDN\r\ and gave some preliminary
results on this concept. Here, we want to examine this notion more
systematically, in particular, we want to
investigate the corresponding convolution. We describe this convolution
both from a combinatorial point of view - by showing its connection with
the lattice of non-crossing partitions - and from an analytic point of
view - by presenting the basic formula for its $R$-transform, which is the
replacement of the classical Fourier-transform. We calculate explicitly
the distributions of the corresponding Gaussian and Poisson law by
a careful examination of the structure of the non-crossing partitions.
\par
Instead of the terms \lq $\psi$-independence' and \lq $\psi$-product'
of \l BSp\r, we will use here the more precise expressions
\lq conditionally free' and
\lq conditionally free product', or just the abbreviation \lq c-free'.
\par
Let us start with a motivation for our concept of
\lq c-freeness'.
Consider a group $G=\ast_{i\in I}G_i$ which is the free product of
groups $G_i$ ($i\in I$), i.e. each element $g\not=e$ of $G$ can
be written uniquely in the form $g=g_1\dots g_n$, where
$e\not= g_j\in G_{i(j)}$ and $i(1)\not=i(2)\not=\dots\not= i(n)$.
To see the nature of this
decomposition of $G$ more clearly, we state it in a more abstract way by using
the functions $\psi_i=\delta_e$ on $G_i$, i.e.
$\psi_i:G_i\to\CC$ with ($g\in G_i$)
 $$\psi_i(g)=\cases 1,& g=e\\
0,& g\not=e.\endcases$$
Then the above decomposition has
the form: Each element $g\not= e$ can be written as $g=g_1\dots g_n$,
where $g_j\in G_{i(j)}$, $i(1)\not=i(2)\not=\dots\not=i(n)$ and
$\psi_{i(j)}(g_j)=0$ for all $j=1,\dots,n$.\par
If we are now given functions $\ff_i:G_i\to\CC$ with $\ff_i(e)=1$, then
we can form their c-free product in the canonical way, namely we
define a new function
$\ff=\ast_{i\in I}(\ff_i,\psi_i)$ on $G$
by $\ff(e)=1$ and
 $$\ff(g):=\ff_{i(1)}(g_1)\dots \ff_{i(n)}(g_n),$$
if $g\not=e$ has the above representation.
The key property of this construction is the fact that, if the
$\ff_i$ are positive definite on $G_i$, then $\ff$ is positive definite
on $G=\ast_{i \in I}G_i$, see \l Boz1,Boz2\r.\par
As an example of such a c-free product one can
take each $G_i$ as a copy of $\ZZ$ and $\ff_i$ as
$\ff_i(k)=\exp(-\lambda \vert k\vert)$ ($k\in\ZZ$) for some $\lambda>0$.
Then $G$ is the free
group on $\vert I\vert$ generators and $\ff$ is given by
$\ff(g):=\exp(-\lambda\vert g\vert)$, where
$g\mapsto \vert g\vert$ is the canonical
length function  on the free group. Since the $\ff_i$ are positive
definite functions on $\ZZ$, this $\ff$ is also positive definite.
This property of the length function on the free group was proven
by Haagerup \l Haa\r.
\par
If we translate the above description of
$\ff$ from groups to group algebras,
then it reads in the following way: Let
$\cA_i:=\CC G_i$
and $\cA:=\CC G$ be the group algebras of $G_i$ and $G$,
respectively. Then,
given linear functionals $\ff_i$ on $\cA_i$ with $\ff_i(1)=1$,
we can define a linear functional $\ff=\ast_{i\in I}(\ff_i,\psi_i)$
on $\cA$ by $\ff(1)=1$ and the
characterizing property
 $$\ff(a_1\dots a_n)=\ff_{i(1)}(a_1)\dots \ff_{i(n)}(a_n),$$
whenever $a_j\in \cA_{i(j)}$, $i(1)\not=i(2)\not=\dots\not= i(n)$
and $\psi_{i(j)}(a_j)=0$,
where $\psi_i$ is now the linear extension of $\delta_e$ to $\cA_i$.
\par
In this formulation it is unnatural to restrict to $\cA_i=\CC G_i$ and
to $\psi_i=\delta_e$,
one can now consider the
above c-free product for arbitrary unital algebras $\cA_i$ and
arbitrary states
$\psi_i$ on $\cA_i$. One of the
main results in \l BSp\r\  was that also in this general case $\ff$ is a
state if the $\ff_i$ are. This was proved by an explicit construction
of the corresponding
c-free product. We will give in Sect. 2 another,
purely algebraic,
proof of this basic fact.
\par
After this basic considerations
we will then switch to the corresponding notion of
c-free convolution, the
main topic of our investigations. Since compactly supported measures $\mu$
on $\RR$ are determined by their moments, such measures can
be identified with states on the polynomial algebra $\CC\la X\ra$. Thus
we can characterize our convolution in the following way. Given pairs of
compactly supported probability measures $(\mu_1,\nu_1)$ and
$(\mu_2,\nu_2)$, we define their c-free convolution
by the following prescription: Consider $\cA_1=\CC\la X_1\ra$ and
$\cA_2=\CC\la X_2\ra$. Then $\cA=\cA_1\ast\cA_2=\CC\la X_1,X_2\ra$.
We have on $\cA_i$ the states $\mu_i$ and $\nu_i$, thus
our construction of a c-free product gives a state
 $\ff=(\mu_1,\nu_1)\ast (\mu_2,\nu_2)$ on
$\cA$. If we restrict this state to $\CC\la X\ra$, where $X=X_1+X_2$,
then the distribution of $X$ determines a measure $\mu$, which we call the
c-free convolution of $(\mu_1,\nu_1)$ and $(\mu_2,\nu_2)$, denoted by
\hbox{$\mu=(\mu_1,\nu_1)\splus (\mu_2,\nu_2)$}. The name \lq c-free
convolution' indicates that
$\mu$ is the distribution of the sum of $X_1$ and $X_2$, which are
distributed according to $\mu_1$ and $\mu_2$ and which are c-free.
If $\nu_i=\mu_i$ ($i=1,2$), then our construction reduces to the
free convolution of Voiculescu \l Voi2\r.
\par
To be able to talk about associativity, we should also define a new measure
$\nu$ and it turns out that the natural candidate for this is the
free convolution $\nu_1\splus\nu_2$ of $\nu_1$ and $\nu_2$, thus
 $$(\mu,\nu)=(\mu_1,\nu_1)\splus(\mu_2,\nu_2),$$
where
 $$\mu=(\mu_1,\nu_1)\splus(\mu_2,\nu_2),\qquad
\nu=\nu_1\splus\nu_2.$$
In Sect. 3, we will examine this c-free convolution from a combinatorial
point of view and show that it is, similarly as in the case of the free
convolution \l Spe2\r, determined by the lattice of non-crossing
partitions.
\par
In Sect. 4, we treat the c-free central and Poisson limit theorem by
a careful analysis of the structure of the non-crossing partitions.
We will thereby derive some combinatorial identities for these partitions
which also have some interest of their own.\par
In Sect. 5, we present a systematic machinery for an analytic description
of c-free convolution, namely the generalization of Voiculescu's
$R$-transform \l Voi2\r.
\vskip1cm
\heading
{\bf 2. Definition and positivity of the c-free product}
\endheading
We shall work in the category of unital $*$-algebras and states. By a
state $\ff$ on a unital $*$-algebra $\cA$ we will always mean a linear
functional $\ff:\cA\to\CC$, which is normalized ($\ff(1)=1$),
hermitian ($\ff(a^*)=\overline{\ff(a)}$ for all $a\in\cA$) and
positive ($\ff(aa^*)\geq 0$ for all $a\in\cA$).\par
Let now $\cA_i$ ($i\in I$) be unital $*$-algebras equipped with a pair
of states $(\ff_i,\psi_i)$. Then we want to define a new state
$\ff=\ast_{i\in I}(\ff_i,\psi_i)$ on the algebraic free product
$\cA=\ast_{i\in I}\cA_i$ (identification of the units is assumed).
Observing that with the decompositions $\cA_i=\CC1\oplus \cA_i^o$, where
$\cA_i^o:=\ker \psi_i$, one can identify $\cA$ as a vector space with
 $$\CC1\oplus\bigoplus_{n=1}^\infty\thinspace
\bigoplus_{i(1)\not=\dots\not= i(n)}
\cA_{i(1)}^o\otimes\dots\otimes \cA_{i(n)}^o,$$
it is clear that we can define uniquely and consistently a linear functional
$\ff=\ast_{i\in I}(\ff_i,\psi_i)$ on $\cA$ by $\ff(1)=1$ and the
following characterization:
 $$\ff(a_1\dots a_n)=\ff_{i(1)}(a_1)\dots \ff_{i(n)}(a_n),$$
whenever
 $$a_j\in\cA_{i(j)},\quad i(1)\not=i(2)\not=\dots\not=i(n),\quad
\psi_{i(j)}(a_j)=0.$$
Such elements $a_1\dots a_n\in
\cA_{i(1)}^o\otimes\dots\otimes \cA_{i(n)}^o$ will be called
{\it elementary elements} in the following. \par
Of course, the main problem is now to see that $\ff$ is positive.
In \l BSp\r, this was proven by an explicit construction of the GNS
representation of $\cA$ with respect to $\ff$. Here, we want to give
a purely algebraic proof of this fact. For this we need
a lemma about the calculation of $\ff$.
\proclaim{Lemma 2.1}
Consider two elementary elements
 $$y_1=a_1^{(1)}\dots a_n^{(1)} \qquad\text{and}\qquad
   y_2=a_1^{(2)}\dots a_m^{(2)}.$$
1) If $a_1^{(1)}$ and $a_1^{(2)}$ do not belong to the same $\cA_i^o$
then
 $$\ff(y_1^*y_2)=\ff(y_1^*)\ff(y_2).$$
2) Consider $a\in\cA_i$ for some $i\in I$. If $a_1^{(1)}$ and
$a_1^{(2)}$ do not belong to $\cA_i^o$ then
 $$\ff(y_1^*ay_2)=\psi_i(a)\ff(y_1^*y_2)-\psi_i(a)\ff(y_1^*)\ff(y_2)
+\ff_i(a)\ff(y_1^*)\ff(y_2).$$
\endproclaim
\demo{Proof}
1) Clear, since
 $$\ff(y_1^*y_2)=\ff(a_n^{(1)*})\dots\ff(a_1^{(1)*})
\ff(a_1^{(2)})\dots\ff(a_m^{(2)}).$$
2) This follows from
 $$\align
\ff(y_1^*ay_2)&=\ff\bigl(y_1^*(a-\psi_i(a)1)y_2\bigr)+\psi_i(a)\ff(y_1^*y_2)\\
&=\ff(y_1^*)\ff_i\bigl(a-\psi_i(a)1\bigr)\ff(y_2)+\psi_i(a)\ff(y_1^*y_2)\\
&=\ff(y_1^*)\ff_i(a)\ff(y_2)-\ff(y_1^*)\psi_i(a)\ff(y_2)+\psi_i(a)\ff(y_1^*
y_2). \qed \endalign$$
\enddemo
\proclaim{Theorem 2.2}
If $\ff_i$ and $\psi_i$ are states for all $i\in I$, then
$\ff=\ast_{i\in I}(\ff_i,\psi_i)$ is a state, too.
\endproclaim
\demo{Proof}
We will show
 $$\ff(x^*x)\geq \vert\ff(x)\vert^2\qquad\text{for all $x\in\cA$.}$$
We can write each $x\in\cA$ in the form
 $$x=\alpha 1+\sum_k \ak_1\dots\ak_{n(k)},$$
where $\alpha\in\CC$ and
$\ak_1\dots\ak_{n(k)}$ are elementary elements (with $n(k)\geq 1$) for
all $k$.
It suffices to prove the asserted inequality for $x$ without term of the
form $\alpha 1$, i.e. we can assume $x$ to be of the form
 $$x=\sum_k \ak\yk$$
with
 $$\ak:=\ak_1\qquad\text{and}\qquad \yk:=\ak_2\dots \ak_{n(k)}.$$
Our proof will be by induction on the length of $x$ (i.e. the
maximal $n(k)$ in the above representation), and we assume now the
validity of the assertion for elements of a smaller length than $x$,
in particular for the $\yk$ and linear combinations of them.\newline
Put now
 $$x_i:=\sum\Sb \text{$k$ with}\\ \ak\in\cA_i\endSb \ak\yk.$$
Then it suffices to prove the assertion for all $x_i$, because this
implies, by the first part of Lemma 2.1.
 $$\align
\ff(x^*x)&=\sum_{i,j}\ff(x_i^*x_j)\\
&=\sum_{i=j}\ff(x_i^*x_i)+\sum_{i\not= j}\ff(x_i^*x_j)\\
&\geq\sum_{i=j}\vert\ff(x_i)\vert^2+\sum_{i\not= j}\overline{\ff(x_i)}
\ff(x_j)\\
&=\vert\sum_i\ff(x_i)\vert^2\\&=\vert\ff(x)\vert^2.
\endalign$$
So let us consider the case of $x_i=\sum\ak\yk$ with all $\ak\in\cA_i$.
Then the second part of Lemma 2.1. gives
 $$\align
\ff(x_i^*x_i)&=\sum_{k,l}\ff(\yks\aks \al\yl)\\
&=\sum_{k,l}\psi_i(\aks\al)\ff(\yks\yl)-\sum_{k,l}\psi_i(\aks\al)
\ff(\yks)\ff(\yl)\\&\quad
+\sum_{k,l}\ff_i(\aks\al)\ff(\yks)\ff(\yl).\endalign$$
By positivity of $\ff_i$ and $\psi_i$ we can write
 $$\ff_i(\aks\al)=\sum_r\overline{\alpha_r^{(k)}}\alpha_r^{(l)},\qquad
  \psi_i(\aks\al)=\sum_r\overline{\beta_r^{(k)}}\beta_r^{(l)}$$
for some $\alpha_r^{(k)},\beta_r^{(k)}\in\CC$.\newline
By using our induction hypothesis for $\sum_k\beta_r^{(k)}\yk$
this implies
 $$\align
\sum_{k,l}\psi_i(\aks\al)\ff(\yks\yl)&=\sum_{k,l,r}\overline{\beta_r^{(k)}}
\beta_r^{(l)}\ff(\yks\yl)\\
&=\sum_r\ff\bigl\lb (\sum_k
\beta_r^{(k)}\yk)^*(\sum_l\beta_r^{(l)}\yl)\bigr\rb\\
&\geq\sum_r\vert\ff(\sum_k\beta_r^{(k)}\yk)\vert^2\\
&=\sum_{k,l}\psi_i(\aks\al)\ff(\yks)\ff(\yl).\endalign$$
But then, again by our induction hypothesis
 $$\align
\ff(x_i^*x_i)&\geq\sum_{k,l}\ff_i(\aks\al)\ff(\yks)\ff(\yl)\\
&=\sum_r\sum_{k,l}\overline{\alpha_r^{(k)}}\alpha_r^{(l)}\ff(\yks)\ff(\yl)\\
&=\sum_r\ff\bigl\lb(\sum_k\alpha_r^{(k)}\ff(\yk))^*(\sum_l\alpha_r^{(l)}
\ff(\yl))\bigr\rb\\
&\geq\sum_{k,l}\ff_i(\aks)\ff_i(\al)\ff(\yks)\ff(\yl)\\
&=\sum_{k,l}\ff(\aks\yks)\ff(\al\yl)\\
&=\vert\ff(x_i)\vert^2.\qed
\endalign$$
\enddemo
\demo{Remarks}
1) If we have $\ff_i=\psi_i$ for all $i\in I$, then we recover the case
of the free product \l Voi1,VDN\r\ and we obtain an algebraic proof
for the positivity also in this case, thus giving a positive answer to
a question posed in \l Spe2\r.\newline
2) If we want to make our construction associative, then we should
extend also the $\psi_i$ to a new state $\psi$ on $\cA$. It is clear that
$\psi$ should be the free product of the $\psi_i$, in our notations
$\psi:=\ast_{i\in I}(\psi_i,\psi_i)$.
This together with
$\ff:=\ast_{i\in I}(\ff_i,\psi_i)$
will be denoted by
 $$(\ff,\psi)=\ast_{i\in I}(\ff_i,\psi_i)$$
(not to be confused with our notation of a symmetrized product in
\l BSp\r). With these definitions one gets directly the
associativity of our c-free product: If $I=I_1\cupdot I_2$ with
$I_1\cap I_2=\emptyset$, then
 $$\ast_{i\in I}(\ff_i,\psi_i)=\bigl\{\ast_{i\in I_1}(\ff_i,\psi_i)
\bigr\}\ast \bigl\{\ast_{i\in I_2}(\ff_i,\psi_i)\bigr\}.$$
3) Commutativity of our construction is clear.
\newline
4) Cabanal-Duvillard \l CDu\r\ introduced a generalization of our
construction from two to infinitely many states. However, his product
ceases to be associative.
\enddemo
\vskip1cm
\heading
{\bf 3. Combinatorial description of the c-free convolution}
\endheading
Let $\cM$ be the set of all compactly supported probability measures on
$\RR$. Since such a measure $\mu$ is determined by its moments we can
identify it with a state on the $*$-algebra $\CC\la X\ra$ (where $X^*=X$)
via
 $$\mu(X^n)=\int t^n d\mu(t)\qquad (n\geq 0).$$
Let now $\mu_1,\mu_2,\nu_1,\nu_2\in\cM$ be given. We identify $\mu_i$,
$\nu_i$ with states on $\CC\la X_i\ra$ ($i=1,2$) and get, by our results
from Sect. 2, the c-free product $\ff=(\mu_1,\nu_1)\ast (\mu_2,\nu_2)$ on
$\CC\la X_1\ra\ast\CC\la X_2\ra=\CC\la X_1,X_2\ra$ (the latter being the
algebra of polynomials in the non-commuting variables $X_1$ and $X_2$).
The c-free convolution
 $$\mu=(\mu_1,\nu_1)\sqp(\mu_2,\nu_2)\in\cM$$
is then given as the distribution of $X:=X_1+X_2$, i.e.
 $$\int t^nd\mu(t)=\mu(X^n)=\ff((X_1+X_2)^n)\qquad (n\geq 0).$$
For $\mu_i=\nu_i$ ($i=1,2$) this reduces to the free convolution of
Voiculescu \l VDN\r.\par
As in Remark 2 of Sect. 2, we define also a measure $\nu$ as the free
convolution of $\nu_1$ and $\nu_2$, i.e.
 $$\nu=(\nu_1,\nu_1)\sqp(\nu_2,\nu_2)=\nu_1\sqp\nu_2,$$
and denote this situation by
 $$(\mu,\nu)=(\mu_1,\nu_1)\sqp(\mu_2,\nu_2).$$
Then our mapping $\sqp:\cM^2\times\cM^2\to\cM^2$ is commutative and
associative.\par
Our aim is now to extend the combinatorial description of the free
convolution with the help of the lattice of non-crossing partitions
\l Spe2\r\ to our case.\par
For a $\mu\in\cM$ we denote its moments by
 $$m_n(\mu):=\mu(X^n)=\int t^n d\mu(t),$$
and we want to understand (at least in principle) the connection between
 $$\bigl(m_n(\mu),m_n(\nu)\bigr)_{n\in\NN}\qquad\text{and}\qquad
\bigl(m_n(\mu_1),m_n(\nu_1)\bigr)_{n\in\NN},
\bigl(m_n(\mu_2),m_n(\nu_2)\bigr)_{n\in\NN}.$$
As in the case of the free convolution this connection is quite
complicated and it is advantageous to introduce new quantities, called
cumulants, which linearize the convolution. These cumulants are connected
with the notion of non-crossing partitions.
\demo{Definition}
Let $\cV=\{V_1,\dots,V_p\}$ be a partition of the linear ordered set
\linebreak
$\{1,\dots,n\}$, i.e. the $V_i\not=\emptyset$ are ordered and disjoint
sets whose union is $\{1,\dots,n\}$. Then $\cV$ is called {\it
non-crossing} if $a,c\in V_i$ and $b,d\in V_j$ with $a<b<c<d$ implies
$i=j$.\newline
The sets $V_i\in\cV$ are called {\it blocks}.
A block $V_i$ of a non-crossing partition $\cV=\{V_1,\dots,V_p\}$ is
called {\it inner}, if there exists a $V_j\in\cV$ and $a,b\in V_j$ such
that
 $a<v<b$ for at least one (and hence for all) $v\in V_i$. A block $V_i\in\cV$
which is not inner is called {\it outer}.
\newline
We will denote the set of all non-crossing partitions of the set
$\{1,\dots,n\}$ by $NC(n)$. By $NC_2(2n)$ we denote those non-crossing
partitions $\cV=\{V_1,\dots,V_n\}\in NC(2n)$ where each block $V_i\in\cV$
consists of exactly two elements.
\enddemo
The notion of non-crossing partition was introduced by Kreweras \l Kre\r,
the distinction between outer and inner blocks was considered in \l BSp\r.
\par
After this preparations we can now introduce the notion of cumulants.
For the description of \hbox{$\nu=\nu_1\sqp\nu_2$} we have to use the
free or non-crossing  cumulants $r_n=r_n(\nu)$,
defined recursively in terms of the
moments $m_n=m_n(\nu)$ by \l Spe2,NSp1,NSp2\r
 $$m_n=\sum_{k=1}^n\thinspace\sum\Sb l(1),\dots,l(k)\geq 0\\
l(1)+\dots+l(k)=n-k\endSb r_k m_{l(1)}\dots m_{l(k)}.$$
This definition may be indicated symbolically by
 $$\boxed{\text{diagram 1}}$$
and it is equivalent to
 $$m_n=\sum\Sb \cV=\{V_1,\dots,V_p\}\\ \in NC(n)\endSb r_{V_1}\dots
r_{V_p}=\sum_{\cV\in NC(n)} \prod_{V_l\in\cV} r_{V_l}\tag 1$$
or
 $$r_n=\sum\Sb \cV=\{V_1,\dots,V_p\}\\ \in NC(n)\endSb m_{V_1}\dots
m_{V_p}\cdot\mu(\cV,1_n)=\sum_{\cV\in NC(n)}\mu(\cV,1_n)
\prod_{V_l\in\cV} m_{V_l},\tag 2$$
where we
have used the notation $m_V:=m_{\vert V\vert}$ and $r_V:=r_{\vert V\vert}$
for some set $V$ (with $\vert V\vert$ being the number of elements in $V$).
The function $\mu(\cV,1_n)$ is the M\"obius function of the lattice
of non-crossing partitions and is just determined by resolving $(1)$
for the $r_n$ in terms of the $m_i$.
\par
Free convolution is then described by \l Voi2,Spe2\r
 $$r_n(\nu_1\sqp\nu_2)=r_n(\nu_1)+r_n(\nu_2)\qquad\text{for all $n\geq 1$.}$$
For our c-free convolution we have, for a given pair $(\mu,\nu)$, to
introduce, in addition to $r_n=r_n(\nu)$, also c-free cumulants
$R_n=R_n(\mu,\nu)$, which do not only depend on the moments of $\mu$ but
also on those of $\nu$. The most instructive definition is again by
recursion, namely
 $$
m_n(\mu)=\sum_{k=1}^n\thinspace \sum\Sb l(1),\dots,l(k)\geq 0\\
l(1)+\dots+l(k)
=n-k\endSb R_k\cdot m_{l(1)}(\nu)
\cdot\dots\cdot m_{l(k-1)}(\nu)\cdot m_{l(k)}(\mu),
 $$
pictorially\newline
 $$\boxed{\text{diagram 2}}$$
Note that the \lq inner' moments are given by $\nu$, only the \lq outer'
one is connected with $\mu$.
Of course, the free cumulants are recovered from this by
$r_n(\nu)=R_n(\nu,\nu)$.\par
The above definition is equivalent to a generalization of $(1)$, namely
 $$m_n(\mu)=\sum_{\cV\in NC(n)} \prod\Sb V_l\in\cV\\ \text{$V_l$ inner}
\endSb r_{V_l} \prod\Sb V_k\in\cV\\ \text{$V_k$ outer}\endSb
R_{V_k}.$$
The following example shows that
the analogue of formula $(2)$ is not true for the c-free cumulants.
\demo{Example}
We have
 $$\align
m_3(\nu)&=r_3(\nu)+2r_2(\nu)\cdot r_1(\nu)+r_1(\nu)^3+r_2(\nu)\cdot r_1(\nu)
  \\
m_3(\mu)&=R_3(\mu,\nu)+2R_2(\mu,\nu)\cdot R_1(\mu,\nu)+R_1(\mu,\nu)^3+
R_2(\mu,\nu)\cdot r_1(\nu),\endalign$$
but
 $$\align
r_3(\nu)&=m_3(\nu)-2m_2(\nu)\cdot m_1(\nu)-m_2(\nu)\cdot m_1(\nu)+
2m_1(\nu)^3\\
R_3(\mu,\nu)&=m_3(\mu)-2m_2(\mu)\cdot
m_1(\mu)-m_2(\mu)\cdot m_1(\nu)+m_1(\mu)^3+
m_1(\mu)^2\cdot m_1(\nu).\endalign$$
\enddemo
But nevertheless we have the following theorem.
\proclaim{Theorem 3.1}
The c-free convolution
 $$(\mu,\nu)=(\mu_1,\nu_1)\sqp(\mu_2,\nu_2)$$
is described by
 $$r_n(\nu)=r_n(\nu_1)+r_n(\nu_2)$$
and
 $$R_n(\mu,\nu)=R_n(\mu_1,\nu_1)+R_n(\mu_2,\nu_2)$$
for all $n\geq 1$.
\endproclaim
\demo{Proof}
The proof follows the same line of argueing as in \l Spe2\r. Given
$(\ff,\psi)$ on some unital $*$-algebra
$\cA$, we define more general cumulant
functions $r=(r_n)$ and $R=(R_n)$ with
 $$r_n,R_n:\undersetbrace \text{$n$-times} \to {\cA\times\dots\times\cA}
\to \CC \qquad (n\geq 1)$$
by
 $$\multline
\psi(a_1\dots a_n)=\sum_{k=0}^{n-1}\thinspace\sum_{1<l(1)<\dots<l(k)\leq n}
r_{k+1}\lb a_1,a_{l(1)},\dots,a_{l(k)}\rb\\
\psi(a_2\dots a_{l(1)-1})\dots\psi(a_{l(k-1)+1}\dots a_{l(k)-1})
\psi(a_{l(k)+1}\dots a_n)\endmultline\tag A$$
and
 $$\multline
\ff(a_1\dots a_n)=\sum_{k=0}^{n-1}\thinspace\sum_{1<l(1)<\dots<l(k)\leq n}
R_{k+1}\lb a_1,a_{l(1)},\dots,a_{l(k)}\rb\\
\psi(a_2\dots a_{l(1)-1})\dots\psi(a_{l(k-1)+1}\dots a_{l(k)-1})
\ff(a_{l(k)+1}\dots a_n)\endmultline\tag B$$
for all $a_1,\dots,a_n\in\cA$. These equations can recursively be resolved
for the definition of $r_n\lb a_1,\dots,a_n\rb$ and $R_n\lb a_1,\dots,a_n\rb$.
\newline
Let now $(\mu_i,\nu_i)$ on $\cA_i=\CC\la X_i\ra$ for $i=1,2$ be given. Then
we obtain in the above way the functions $r(\mu_i)$ and $R(\mu_i,\nu_i)$
on $\bigcup_{n=1}^\infty \cA_i^{\times n}$.
On $\bigcup_{n=1}^\infty (\cA_1\cup \cA_2)^{\times n}\subset
\bigcup_{n=1}^\infty (\cA_1\ast \cA_2)^{\times n}$ we define their
direct sum
 $$r:=r(\mu_1)\oplus r(\mu_2)\qquad\text{and}\qquad R:=R(\mu_1,\nu_1)
\oplus R(\mu_2,\nu_2)$$
by
 $$r_n\lb a_1,\dots,a_n\rb=\cases
r_n(\mu_1)\lb a_1,\dots,a_n\rb,&\text{if all $a_i\in\cA_1$}\\
r_n(\mu_2)\lb a_1,\dots,a_n\rb,&\text{if all $a_i\in\cA_2$}\\
0,&\text{otherwise}\endcases$$
and
 $$R_n\lb a_1,\dots,a_n\rb=\cases
R_n(\mu_1,\nu_1)\lb a_1,\dots,a_n\rb,&\text{if all $a_i\in\cA_1$}\\
R_n(\mu_2,\nu_2)\lb a_1,\dots,a_n\rb,&\text{if all $a_i\in\cA_2$}\\
0,&\text{otherwise}\endcases$$
for all $a_1,\dots,a_n\in\cA_1\cup\cA_2\subset\cA_1\ast\cA_2$.
Note that there is no ambiguity in this definition because in the
case that some $a_i\in\cA_1\cap \cA_2=\CC 1$, both values, $R_n(\mu_1,\nu_1)$
and $R_n(\mu_2,\nu_2)$, are the same.
\newline
Now we use the recursion formulas $(A)$ and $(B)$ for the definition of
the states $\psi$ and $\ff$ on $\cA=\cA_1\ast\cA_2$. One has to check
that this is well-defined because there are different possibilities for
writing elements $a\in\cA$ as sums of products $a_1\dots a_n$ with
$a_1,\dots,a_n\in\cA_1\cup\cA_2$. But since this ambiguity comes only
from relations inside $\cA_1$ and relations inside $\cA_2$, which are
respected by $r$ and $R$ (because they are respected by $r(\nu_i)$ and
$R(\mu_i,\nu_i)$), no problem occurs; for more details on this, see
\l Spe2\r.\newline
It only remains to see that $(\ff,\psi)$ on $\cA=\CC\la X_1,X_2\ra$ is
indeed the c-free product of $(\mu_1,\nu_1)$ and $(\mu_2,\nu_2)$, i.e. we
have to check that it
fulfills the characterizing property of the c-free product. For $\psi$
this follows from the results of \l Spe2\r.
So consider $a\in\cA$ of the form
$a=a_1\dots a_n$ with $a_j\in\cA_{i(j)}$,
 $$
i(1)\not=i(2)\not=\dots\not=i(n),\quad \nu_{i(j)}(a_j)=
0.$$
Note that in $(B)$, because of the definition of $R$ and the fact that
$\psi$ is the free product of $\psi_1$ and $\psi_2$, only the term with
$k=0$ survives, i.e.
 $$\ff(a_1\dots a_n)=R_1\lb a_1\rb\cdot \ff(a_2\dots a_n)=
\ff_{i(1)}(a_1)\cdot \ff(a_2\dots a_n),$$
which gives, by induction, the wanted factorization for $\ff$.\newline
To get the assertion of the theorem, one has now to use the definition
of $R$ as the direct sum of $R(\mu_1,\nu_1)$ and $R(\mu_2,\nu_2)$
 $$\align
R_n(\mu,\nu)&\hateq R_n\lb X_1+X_2,\dots,X_1+X_2\rb\\
&=R_n(\mu_1,\nu_1)\lb X_1,\dots,X_1\rb+R_n(\mu_2,\nu_2)\lb X_2,\dots,X_2\rb\\
&\hateq R_n(\mu_1,\nu_1)+R_n(\mu_2,\nu_2),\endalign
 $$
and the same for $r$. \qed
\enddemo
\demo{Remarks}
1) The description of the c-free convolution in terms of cumulants can,
analogously to the free case \l Spe2\r, be generalized to a description
of the c-free product. Indeed, in our proof we had to use the corresponding
machinery for the c-free product on $\CC\la X_1,X_2\ra=\CC\la X_1\ra\ast
\CC\la X_2\ra$.\newline
2) An interesting special case of the c-free convolution is given if we
put $\nu_i=\delta_0$. Then only outer sets survive in the definition of
the c-free cumulants. This leads to a description in terms of interval
partitions, which were introduced by von Waldenfels \l vWa\r.
The corresponding convolution \hbox{$(\mu,\delta_0)=(\mu_1,\delta_0)\splus
(\mu_2,\delta_0)$} shares a lot of properties with the usual and the
free convolution. This \lq boolean' convolution was investigated in
\l Wor\r, the results will be published in \l SpW\r.
\enddemo
\vskip1cm
\heading
{\bf 4. Limit theorems for the c-free convolution}
\endheading
To become familiar with the connection between non-crossing partitions and
the c-free convolution, we will now calculate quite
explicitly the c-free central and Poisson limit distribution. A more
systematic machinery for the treatment of such questions will be
presented in the next section.\par
We will see (comp. \l BSp\r) that the moments of the limit distributions
are calculated with the help of the partitions in $NC_2(2n)$ or $NC(n)$.
Thus, before presenting the limit theorems, we collect all relevant
information on the combinatorics of the respective partitions in two
lemmas. These combinatorial statements have also some interest of
their own. Although there has been an increasing interest in the lattice
of non-crossing partitions
in the last time \l Ede1,Ede2,Pou,Sim,SiU,Bia,Nic\r,
we have not found any investigation on this subject related to the distinction
between \lq inner' and \lq outer'.\par
First, for the central limit theorem, we have to consider $\NCzn$. We will
need the numbers ($n\in\NN$, $0\leq k\leq n$)
 $$\align c^n&:=\#\NCzn\\
a^n_k&:=\#\{\cV\in\NCzn\mid \text{the number of inner sets of $\cV$ is equal
to $k$}\}.\endalign$$
Of course, we have $a^n_n=0$.\par
For the investigation of these quantities it is advantageous to use the
well known bijection between partitions $\cV\in\NCzn$ and $n$-Catalan
paths $\Lambda$ (see, e.g., \l HiP\r).
An $n$-Catalan path $\Lambda\hateq\{s_1,\dots,s_{2n}\}$
is a graph in $\ZZ^2$, starting at $(0,0)$, ending at $(n,n)$, with
possible steps $s_i=(0,1)$ or $s_i=(1,0)$ ($i=1,\dots,2n$), such that the
path does not lie above the diagonal. The above bijection
is given as follows: To each $\cV=\{V_1,\dots,V_n\}\in\NCzn$ we assign a
$\Lambda(\cV)\hateq\{s_1,\dots,s_{2n}\}$ in the way that $s_i=(1,0)$ if
$i$ is the first element in one of the $V_j$, and $s_i=(0,1)$ if $i$
is the second element in one of the $V_j$. The number of outer sets of
$\cV$ corresponds thereby to the number of points $(i,i)$ ($1\leq i\leq n$),
where $\Lambda(\cV)$ meets the diagonal.
\demo{Example}
For
 $\cV=\{(1,4),(2,3),(5,6)\}$ we have
 $$\Lambda(\cV)=\{(1,0),(1,0),(0,1),(0,1),(1,0),(0,1)\},$$ which
corresponds to the following graph:
  $$\boxed{\text{diagram 3}}$$
\enddemo
It is a well known fact \l HiP\r\ that the number of all $n$-Catalan paths
is given by the
 $$\text{Catalan number}\qquad c_n:=\frac 1n\binom {2n}{n-1},$$
hence
 $$c^n=c_n=\frac 1n\binom{2n}{n-1}=\frac 1{n+1}\binom{2n}n.$$
This follows quite easily from the recursion formula
 $$c^n=\sum_{k=1}^n c^{k-1}c^{n-k},\qquad\text{where $c^0:=1$,}$$
which is the recursion for the Catalan numbers.\par
It seems that $a_k^n$ has not received any interest so far. We collect
their basic properties in the next lemma.
\proclaim{Lemma 4.1}
i) We have for $n\geq 1$
  $$a_0^n=1.$$
ii) We have for $n\geq 2$
 $$a_{n-1}^n=a_{n-2}^n=c^{n-1}=\frac 1n\binom{2n-2}{n-1}.$$
iii) We have for $n\geq2$ and $0\leq k\leq n-2$
 $$a_k^n+a_{k+1}^{n-1}=a_{k+1}^n.$$
\endproclaim
\demo{Proof}
i) There is only one $n$-Catalan path which meets the diagonal $n$-times,
namely $\Lambda\hateq\{(1,0),(0,1),(1,0),(0,1),\dots\}$.\newline
ii) Shifting the diagonal one unit to the right induces a bijection
between the set of all $n$-Catalan paths which meet the diagonal once and the
set of all $(n-1)$-Catalan paths. Hence $a_{n-1}^n=c^{n-1}$.\newline
For $a_{n-2}^n$ we have, denoting with $(k,k)$ the first intersection point
with the diagonal,
 $$a_{n-2}^n=\sum_{k=1}^{n-1}a_{k-1}^k a_{n-k-1}^{n-k}=\sum_{k=1}^{n-1}
c^{k-1}c^{n-k-1}=c^{n-1},$$
by the recursion formula for the Catalan numbers.\newline
iii) We prove this by induction on $n$. For $n=2$, the assertion is
true, namely for $k=0$ we have
 $$a_0^2+a_1^1=1+0=1=a_1^2.$$
Now assume the assertion to be true for all $n'$ with $2\leq n'<n$. We want
to show it for $n$.\newline
First, consider $k$ with $0\leq k\leq n-4$. Again we use the general
decomposition
 $$a_{k+1}^n=\sum_{l=1}^{k+2}a_{l-1}^l a_{k-l+2}^{n-l},\tag *$$
which results from the splitting of an $n$-Catalan path into two parts,
the first one from $(0,0)$ to its first intersection point $(l,l)$ with
the diagonal (this part thus gives rise to $l-1$ inner sets) and the
remaining $(n-l)$-Catalan path, which has to produce the remaining
$(k+1)-(l-1)$ inner sets. The decomposition $(*)$ is true for all
$n\geq 2$ and $k$ with $0\leq k\leq n-3$. Since $0\leq k\leq n-4$, we have
$0\leq k-l+1\leq n-l-3\leq n-l-2$ and, for all $l$ with $1\leq l\leq k+1$,
we can use our induction hypothesis for $n'=n-l$ to obtain
 $$\align
a_{k+1}^n&
=\sum_{l=1}^{k+1} a_{l-1}^l a_{(k-l+1)+1}^{n-l}+a_{k+1}^{k+2}\cdot 1\\
&=\sum_{l=1}^{k+1}a_{l-1}^l(a_{k-l+1}^{n-l}+a_{k-l+2}^{n-l-1})+
a_{k+1}^{k+2}\cdot 1\\
&=\sum_{l=1}^{k+1}a_{l-1}^l a_{k-l+1}^{n-l}+\sum_{l=1}^{k+2}
a_{l-1}^l a_{k-l+2}^{n-l-1}\\
&=a_k^n+a_{k+1}^{n-1},\endalign$$
the last equality again by application of $(*)$.\newline
Now consider $k=n-3$. Then the
same arguments as before apply, but now $a_{k-l+2}^{n-l-1}=0$ for
$l=1,\dots,k+1$, thus
 $$a_{k+1}^n=\sum_{l=1}^{k+1}a_{l-1}^l a_{k-l+1}^{n-l}+a_{k+1}^{k+2}=
a_k^n+a_{k+1}^{n-1}.$$
For $k=n-2$, the assertion reduces to ii), because $a_{k+1}^{n-1}=
a_{n-1}^{n-1}=0$. \qed
\enddemo
For the treatment of the c-free Poisson distribution we will need some
specific information on the combinatorics of the sets $NC(n)$, namely we
will use ($n\geq 1$, $1\leq k\leq n$, $0\leq l\leq n-1$)
 $$\align
t_k^n&:=\#\{\cV\in NC(n)\mid\text{$\cV$ consists precisely of $k$ sets}\}\\
s^n_{k,l}&:=\#\{\cV\in NC(n)\mid \text{$\cV$ consists precisely of $k$
outer and $l$ inner sets}\}.\endalign$$
In addition, we define $t^n_0:=0$ for $n\geq 1$ and $t_0^0:=1$.
Similarly, we put $s_{k,l}^n:=0$ if the indices are out of their natural
domain, with the only exception $s_{0,0}^0:=1$.
\proclaim{Lemma 4.2}
i) We have for $n\geq1$ and $1\leq k\leq n$
 $$t_k^n=t_{k-1}^{n-1}+\sum_{r=2}^n\sum_{i=1}^{r-1}t_i^{r-1}t_{k-i}^{n-r}.$$
ii) We have for $n\geq 2$ and $0\leq l\leq n-1$
 $$s_{1,l}^n=t^{n-1}_{l+1}.$$
iii) We have for $n\geq 1$ and $k,l\geq 0$
 $$s_{k+1,l}^n=\sum_{r=1}^n\sum_{j=0}^l s_{1,j}^r s_{k,l-j}^{n-r}.$$
\endproclaim
\demo{Proof}
i) Let $\cV=\{V_1,\dots,V_k\}\in NC(n)$ consist of $k$ sets, with
$1\in V_1$. Then there are two disjoint possibilities: either $V_1=(1)$
or $V_1\not=(1)$. In the first case, $\cV\mapsto \cV\backslash (1)$ gives
a bijection onto all non-crossing partitions of $\{2,\dots,n\}$ consisting
of $k-1$ sets. In the second case, let $r\not=1$ be the maximal element
of $V_1$. Then, removing $r$ from $V_1$, $\cV$ splits into a non-crossing
partition of $\{1,\dots,r-1\}$ (which may consist of $i$ sets, where
possibly $1\leq i\leq r-1$) and a non-crossing partition of
$\{r+1,\dots,n\}$ (which has to consist of the remaining $k-i$ sets).
If $r=n$, then $k=i$, and we need the special definition $t_0^0:=1$. The
formula is also true for $n=1$, since then $t_1^1=t_0^0=1$.\newline
ii) If $\cV\in NC(n)$ has $l$ inner sets
and only one outer set $V_1$, then $1,n\in V_1$ and the
removing of $n$ ($n\not=1$) gives a bijection onto all non-crossing
partitions of $\{1,\dots,n-1\}$ consisting of $l+1$ sets.\newline
iii) Let $r$ be the maximal element of the first set $V_1$ in
$\cV\in NC(n)$. Then $\cV$ decomposes into a non-crossing partition
of $\{1,\dots,r\}$ with $V_1$ as the only outer set (and possibly $j$
inner sets) and a non-crossing partition of $\{r+1,\dots,n\}$ which has
to yield the remaining $k$ outer and $l-j$ inner sets. If $k=0$ and
$r=n$, then $j=0$, and we need $s_{0,0}^0=1$.
\qed
\enddemo
\demo{Remark}
Kreweras \l Kre\r\ gives the following explicit formula for $t_k^n$
 $$t_k^n=\frac{(n-1)!n!}{(k-1)!k!(n-k)!(n-k+1)!},$$
but for our investigations the recurrence formula of our lemma is much
more useful.
\enddemo
Now we have finished the presentation of all needed combinatorial tools
and we can start our investigations on limit theorems for the c-free
convolution.\par
Let us denote, for $\lambda>0$,
by $D_\lambda$ the dilation of probability measures on
$\RR$ by the factor $\lambda$, i.e.
 $$(D_\lambda\mu)(A):=\mu(\lambda^{-1}A)\qquad\text{for $A\subset\RR$
measurable,}$$
and
 $$D_\lambda(\mu,\nu):=(D_\lambda \mu,D_\lambda\nu).$$
Under the weak convergence
 $$\wlim_{N\to\infty}(\mu_N,\nu_N)=(\mu,\nu)$$
we will understand the componentwise weak convergence
 $$\wlim_{N\to\infty}\mu_N=\mu\qquad\text{and}\qquad
   \wlim_{N\to\infty}\nu_N=\nu.$$
\proclaim{Theorem 4.3 (c-free central limit theorem)}
Let $(\mu,\nu)\in\cM^2$ with
 $$\mu(X)=\nu(X)=0\qquad\text{and}\qquad
\mu(X^2)=\alpha^2,\quad \nu(X^2)=\beta^2\quad (\alpha,\beta\geq 0)$$
be given. Then we have
 $$\wlim_{N\to\infty} D_{\sqrt{1/N}}\bigl\{
\undersetbrace \text{$N$-times}\to{(\mu,\nu)\splus\dots\splus(\mu,\nu)}
\bigr\}=(\nu_{\alpha,\beta},\nu_{\beta,\beta}),$$
where
 $$\nu_{\alpha,\beta}=c(\delta_{\alpha^2/\sqrt{\alpha^2-
\beta^2}}+
\delta_{-\alpha^2/\sqrt{\alpha^2-
\beta^2}})+\tilde \nu_{\alpha,\beta},$$
with  $$\align
   c&=\cases \frac 14\frac{\alpha^2-2\beta^2}{\alpha^2
-\beta^2},& 0\leq\frac {\beta^2}{\alpha^2}\leq\frac 12\\
0,&
\frac 12\leq\frac {\beta^2}{\alpha^2}\endcases\\
d\tilde\nu_{\alpha,\beta}(t)&=\chi_{\lb -2\beta,2\beta\rb}(t)
\frac 1{2\pi} \frac{\alpha^2\sqrt{4\beta^2-t^2}}{\alpha^4-(\alpha^2-
\beta^2)t^2}dt.\endalign$$
In particular
 $$d\nu_{\beta,\beta}(t)=\chi_{\lb-2\beta,2\beta\rb}(t)
                 \frac 1{2\pi\beta^2}\sqrt{4\beta^2-t^2}dt.$$
\endproclaim
\demo{Remark}
Of course, the statement about the convergence of the second
component is nothing else but the free
central limit theorem \l Voi2,VDN,Spe1,Maa,Gir\r.
\enddemo
\demo{Proof}
Since $\nu_{\alpha,\beta}$ and $\nu_{\beta,\beta}$ have compact support, it
suffices to check that the moments of
 $D_{\sqrt{1/N}}\{(\mu,\nu)\splus\dots\splus(\mu,\nu)\}$
converge to the corresponding
moments of $(\nu_{\alpha,\beta},\nu_{\beta,\beta})$.
Note that
 $$r_n(D_\lambda\nu)=\lambda^n r_n(\nu)\qquad\text{and}\qquad
R_n(D_\lambda\mu,D_\lambda\nu)=\lambda^n R_n(\mu,\nu)$$
for all $n\geq 0$. This shows that the limiting measures $(\mh,\nh)$
are determined by
 $$\align r_n(\nh)&=\cases 0,& n\not=2\\
r_2(\nu)=\beta^2,& n=2\endcases\\
R_n(\mh,\nh)&=\cases 0,& n\not=2\\
R_2(\mu,\nu)=\alpha^2,& n=2,\endcases\endalign$$
or in terms of their moments
 $$\align
m_l(\nh)&=\cases 0,&\text{$l$ odd}\\
c^n\beta^{2n},& l=2n\endcases\\
m_l(\mh)&=\cases 0,& \text{$l$ odd}\\
\sum_{k=0}^{n-1} a^n_k \alpha^{2(n-k)}\beta^{2k},& l=2n.\endcases\endalign$$
This formula for the moments of $\mh$ was also derived in
\l BSp\r.\newline
Now consider the generating power series
 $$f(z):=\sum_{n=0}^\infty m_{2n}(\nh)z^{2n},\qquad
F(z):=\sum_{n=0}^\infty m_{2n}(\mh)z^{2n}.$$
The recursion formula for the Catalan numbers yields \l Spe2,VDN\r
 $$\beta^2z^2f(z)^2=f(z)-1,\qquad\text{thus}\qquad
f(z)=\frac{1-\sqrt{1-4\beta^2z^2}}{2\beta^2z^2}.$$
For the determination of $F(z)$, we use part iii) of Lemma 4.1 to
observe
 $$\align
(\alpha^2&-\beta^2)m_{2(n+1)}(\mh)=
\sum_{k=0}^n a^{n+1}_k \alpha^{2(n+2-k)}\beta^{2k}
-\sum_{k=0}^n a^{n+1}_k \alpha^{2(n+1-k)}\beta^{2(k+1)}\\
&=\sum_{k=0}^{n-1} (a_{k+1}^{n+1}-a_k^{n+1})\alpha^{2(n+1-k)}\beta^{2(k+1)}
+a_0^{n+1}\alpha^{2(n+2)}-a_n^{n+1}\alpha^2\beta^{2(n+1)}\\
&=\sum_{k=0}^{n-1} a_{k+1}^{n}\alpha^{2(n+1-k)}\beta^{2(k+1)}+
\alpha^{2(n+2)}
-a_n^{n+1}\alpha^2\beta^{2(n+1)}\\
&=\sum_{k=0}^n a^n_k \alpha^{2(n+2-k)}\beta^{2k}-\alpha^2\beta^{2(n+1)}
a^{n+1}_n\\
&=\alpha^4 m_{2n}(\mh)-\alpha^2 \beta^{2(n+1)} c^n,\endalign$$
which implies
 $$\align
(\alpha^2-\beta^2)F(z)&=
\amb+\sum_{n=0}^\infty m_{2(n+1)}(\mh)z^{2(n+1)}\\
&=\amb+\alpha^4z^2\sum_{n=0}^\infty m_{2n}(\mh)z^{2n}
-\alpha^2\beta^2z^2\sum_{n=0}^\infty c^n (\beta z)^{2n}\\
&=\amb+\alpha^4z^2F(z)-\alpha^2\beta^2z^2 f(z).\endalign$$
This can be resolved for $F(z)$,
 $$
F(z)=\frac{\amb-\alpha^2\beta^2z^2f(z)}{\amb-\alpha^4z^2}
=\frac{\amb-\frac12 \alpha^2(1-\sqrt{1-4\beta^2z^2})}{\amb-\alpha^4
z^2}.$$
In terms of the Cauchy-transform $G(z)$ of $\mh$ this reads
 $$G(z)=\frac 1z F(\frac 1z)=
\frac{z(\frac12 \alpha^2-\beta^2)
+\frac 12\alpha^2 \sqrt{z^2-4\beta^2}}{z^2\amb -\alpha^4}.
 $$
The Stieltjes inversion formula (see, e.g., \l AGl\r) gives then the
distribution as stated in the theorem. \qed
\enddemo
\demo{Remarks}
1) An instructive way to write the Cauchy-transforms\newline
$g(z)=1/zf(1/z)$ and $G(z)=1/zF(1/z)$
of $\nu_{\beta,\beta}$ and $\nu_{\alpha,\beta}$, respectively,
are the following continued fraction expressions
 $$g(z)=\cfrac
1\\z-\cfrac \beta^2\\z-\cfrac \beta^2\\z-\cfrac \beta^2\\z-\dots\endcfrac,
\qquad
G(z)=\frac 1{z-\alpha^2g(z)}=
\cfrac 1\\z-
\cfrac \alpha^2\\z-\cfrac \beta^2\\z-\cfrac \beta^2\\z-\dots\endcfrac.
         $$
These expansions follow directly from the relations
 $$\align \beta^2z^2f(z)f(z)&=f(z)-1\\
          \alpha^2z^2F(z)f(z)&=F(z)-1.\endalign$$
The second identity can be checked with our explicit form of $f$ and
$F$ or it may be derived directly by the recursion formula
 $$m_{2n}(\mh)=\sum_{k=1}^n\alpha^2 m_{2(k-1)}(\nh)m_{2(n-k)}(\mh).$$
2) The sequence of orthogonal polynomials corresponding to $\nu_{\alpha,
\beta}$ satisfies the following recurrence relations:
 $$\align
p_0(x)&=1\\
p_1(x)&=x\\
p_2(x)&=x^2-\alpha^2\\
p_{n+1}(x)&=xp_n(x)-\beta^2p_{n-1}(x)\qquad (n\geq 2).
\endalign$$
For $\alpha^2=\beta^2=1$ we obtain the Tchebyscheff polynomials of the
second kind, whereas for $\alpha^2=1$ and $\beta^2=1/2$ we get the
Tchebyscheff polynomials of the first kind.
\newline
3) It may be interesting to note that in the limit $\alpha,\beta\to\infty$
under the restriction $\beta/\alpha^2=\text{const}=\gamma$, the
distribution $\nu_{\alpha,\beta}$ converges to the Cauchy distribution $\mu$
with density
 $$d\mu(t)=\frac 1\pi\frac \gamma{1+\gamma^2 t^2}dt.$$
4) In Fig. 1, we have plotted the density of $\nu_{\alpha,\beta}$ for
fixed $\beta=1$ and for six different values of $\alpha$.
\enddemo
 $$\boxed{\text{Fig. 1}}$$
\proclaim{Theorem 4.4 (c-free Poisson limit theorem)}
For $\alpha,\beta\geq 0$ define for all $N\geq 1$
 $$\mu_N:=(1-\frac \alpha N)\delta_0+\frac \alpha N \delta_1\qquad\text{and}
\qquad
\nu_N:=(1-\frac \beta N)\delta_0+\frac \beta N \delta_1.$$
Then we have
 $$\wlim_{N\to\infty}\{\undersetbrace \text{$N$-times} \to
{(\mN,\nN)
\splus\dots\splus(\mN,\nN)}\}=(\pi_{\alpha,\beta},\pi_{\beta,\beta}),$$
where
 $$\pi_{\alpha,\beta}=a\delta_0+b\delta_{z_0}
+\tilde\pi_{\alpha,\beta}$$
with
 $$\align
z_0&=\alpha+\frac\alpha{\alpha-\beta}\\
a&=\cases \frac{1-\beta}{1+\alpha-\beta},& 0\leq\beta\leq 1\\
0,& 1\leq\beta\endcases\\
b&=\cases \frac{\beta z_0-\alpha^2}{z_0(\beta-\alpha)},& \alpha\leq\beta
-\sqrt\beta\quad\text{or}\quad \beta+\sqrt\beta\leq\alpha\\
0,& \beta-\sqrt\beta\leq\alpha\leq\beta+\sqrt\beta\endcases
\endalign$$
and
 $$d\tilde\pi_{\alpha,\beta}(t)=\chi_{\lb(1-\sqrt\beta)^2,(1+\sqrt\beta)^2\rb}
(t)\frac1\pi\frac{\alpha\sqrt{4
\beta-\bigl(t-(1+\beta)\bigr)^2}}{2t\lb t(\beta-\alpha)
+\alpha(1-\beta+\alpha)\rb}dt.$$
In particular
 $$\align\pi_{\beta,\beta}&
=\cases (1-\beta)\delta_0+\tilde\pi_{\beta,\beta},& 0\leq\beta\leq1\\
\tilde\pi_{\beta,\beta},& 1\leq\beta \endcases\\
d\tilde\pi_{\beta,\beta}(t)&=
\chi_{\lb(1-\sqrt\beta)^2,(1+\sqrt\beta)^2\rb}(t)
\frac 1{2\pi t}\sqrt{4\beta-\bigl(t-(1+\beta)\bigr)^2}dt.
 \endalign$$
\endproclaim
\demo{Remark}
Again, the statement about the second
component
reduces to the free Poisson limit
theorem \l Maa,VDN,Spe1,Gir\r.
\enddemo
\demo{Proof}
Again, it is sufficient to check the convergence of all moments. Since
for $n\geq 1$
 $$m_n(\nN)=\frac \beta N\qquad\text{and}\qquad m_n(\mN)=\frac\alpha N,$$
we have
 $$r_n(\nN)=\frac \beta N+O(1/N^2)\qquad\text{and}\qquad
R_n(\mN,\nN)=\frac\alpha N+O(1/N^2),$$
from which it follows that the limiting measures $(\mh,\nh)$ are
determined by
 $$\align
r_n(\nh)&=\beta\qquad\text{for all $n\geq 1$}\\
R_n(\mh,\nh)&=\alpha\qquad\text{for all $n\geq 1$},\endalign$$
or equivalently, for all $n\geq 1$,
 $$\align
m_n(\nh)&=\sum_{k=1}^n t^n_k\beta^k\\
m_n(\mh)&=\sum_{k,l\geq 0} s^n_{k,l}\alpha^k\beta^l.\endalign$$
For $\nh$, this gives the free Poisson distribution, see
\l VDN,Maa,Spe1\r.
The formula for the moments of $\mh$ was also derived in
\l BSp\r.
\newline
As before, we want to calculate the generating power series in the
moments
 $$f(z):=\sum_{n=0}^\infty m_n(\nh)z^n\qquad\text{and}\qquad
F(z):=\sum_{n=0}^\infty m_n(\mh)z^n.$$
Since $f(z)$ is of eminent importance for the determination of $F(z)$, we
will briefly derive its form, although this may also be found in
\l Maa,VDN\r. By part i) of Lemma 4.2, we obtain
 $$\align
f(z)&=1+\sum_{n=1}^\infty\bigl(\sum_{k=1}^n t^n_k \beta^k\bigr) z^n\\
&=1+\sum_{n=1}^\infty\sum_{k=1}^n\bigl(t^{n-1}_{k-1}+\sum_{r=2}^n
\sum_{i=1}^{r-1} t_i^{r-1}t_{k-i}^{n-r}\bigr)\beta^kz^n\\
&=1+\beta z\sum_{n=1}^\infty\sum_{k=1}^n t^{n-1}_{k-1}\beta^{k-1}z^{n-1}
+h(z)\\
&=1+\beta z f(z)+h(z),\endalign$$
where
 $$h(z)=z\sum_{n=1}^\infty\sum_{k=1}^n\sum_{r=2}^n\sum_{i=1}^{r-1}
t_i^{r-1}\beta^iz^{r-1}t_{k-i}^{n-r}\beta^{k-i}z^{n-r}
=z(f(z)-1)f(z),$$
thus
 $$f(z)=1+\beta zf(z)+z(f(z)-1)f(z).$$
This can be resolved to give (note $f(0)=1$)
 $$f(z)=\frac {1-(\beta-1)z-\sqrt{\bigl(1-(\beta-1)z\bigr)^2-4z}}{2z}$$
or
 $$g(z)=\frac 1zf(\frac 1z)=
\frac{z+(1-\beta)-\sqrt{\bigl(z-(1+\beta)\bigr)^2-4\beta}}{2z}.$$
For the determination of $F(z)$, we use part iii) and ii) of Lemma 4.2.
We have
 $$\align F(z)&=1+\sum_{n=1}^\infty\sum_{l,k\geq 0} s^n_{k+1,l}
\alpha^{k+1}\beta^l z^n\\
&=1+\sum_{n=1}^\infty\sum_{l,k\geq 0}\sum_{r=1}^n\sum_{j=0}^l
(s^r_{1,j}\alpha \beta^j z^r)(s^{n-r}_{k,l-j}
\alpha^k\beta^{l-j}z^{n-r})\\
&=1+\sum_{j\geq 0}\sum_{r\geq 1}(s^r_{1,j}\alpha\beta^j z^r)F(z)\\
&=1+F(z)\bigl\{\sum_{j\geq 0}s^1_{1,j}\alpha\beta^j z+\frac \alpha\beta z
\sum_{j\geq 0}\sum_{r\geq 2} t^{r-1}_{j+1}\beta^{j+1}z^{r-1}\bigr\}\\
&=1+F(z)\bigl\{\alpha z+\frac\alpha\beta z(f(z)-1)\bigr\},\endalign$$
which implies
 $$F(z)=\frac\beta{\beta-\alpha z(f(z)-1+\beta)}.$$
This yields for the Cauchy-transform $G(z)=1/zF(1/z)$ of $\mh$
after some calculations the expression
 $$G(z)=\frac{z(2\beta-\alpha)+\alpha(1-\beta)-\alpha
\sqrt{\bigl(z-(1+\beta)\bigr)^2
-4\beta}}{2z\lb z(\beta-\alpha)+\alpha(1-\beta+\alpha)\rb}.$$
The Stieltjes inversion formula \l AGl\r\ gives then, after some
computations, the distribution as stated in the theorem. \qed
\enddemo
\demo{Remarks}
1) Again, it is quite instructive to write the Cauchy-transforms as
infinite continued fractions, namely
 $$g(z)=\cfrac 1\\ z+(1-\beta)-\cfrac z\\ z+(1-\beta)-\cfrac z\\
z+(1-\beta)-\cfrac z\\ \dots\endcfrac$$
and
 $$G(z)=\cfrac 1\\z+\frac\alpha\beta(1-\beta)-\cfrac \frac \alpha\beta z\\
z+(1-\beta)-\cfrac z\\z+(1-\beta)-\cfrac z\\\dots\endcfrac.$$
2) Note that our formula for $G(z)$ in \l BSp\r\ was wrong. \newline
3) In Fig. 2, we show the Poisson limit distribution $\pi_{\alpha,
\beta}$ for
$\alpha=1$ and for six different values of $\beta$.
\enddemo
 $$\boxed{\text{Fig. 2}}$$
\vskip1cm
\heading
{\bf 5. Analytic description of the c-free convolution}
\endheading
In Sect. 3, we described the c-free convolution from a combinatorial
point of view by presenting the connection between moments and free and
c-free cumulants; the convolution is then characterized by the fact that
the cumulants are linear under convolution. \par
For an analytic description one wants to translate this connection into
a functional relation between the corresponding power series, i.e.
instead of a collection of moments or cumulants one prefers to deal with
one respective analytic function containing the same information. This
has the advantage that an analytic machinery is usual more powerful
than a mere combinatorial description and it may serve as a starting point
for the treatment of measures with unbounded support.\par
Thus, given a pair $(\mu,\nu)\in\cM^2$, we define the following
power series (formally, we put $r_0=R_0=0$)
 $$\align
A(z)&:=\sum_{n=1}^\infty r_n(\nu)z^n\\
B(z)&:=\sum_{n=0}^\infty m_n(\nu)z^n=1+\sum_{n=1}^\infty m_n(\nu)z^n\\
C(z)&:=\sum_{n=1}^\infty R_n(\mu,\nu)z^n\\
D(z)&:=\sum_{n=0}^\infty m_n(\mu)z^n=1+\sum_{n=1}^\infty m_n(\mu)z^n.
\endalign$$
Since $r_n$ and $R_n$ are additive under c-free convolution, one has
for $(\mu,\nu)=(\mu_1,\nu_1)\splus(\mu_2,\nu_2)$
 $$\align
A_\nu(z)&=A_{\nu_1}(z)+A_{\nu_2}(z)\\
C_{(\mu,\nu)}(z)&=C_{(\mu_1,\nu_1)}(z)+C_{(\mu_2,\nu_2)}(z)\endalign$$
and it remains to derive the connection between $A(z)$ and $C(z)$ on one
side and $B(z)$ and $D(z)$ on the other side. Since $\nu=\nu_1\splus\nu_2$
is nothing else than the free convolution, the relation between
$A(z)$ and $B(z)$ is given in \l Voi2,Spe2\r.
\proclaim{Theorem 5.1}
With the above definitions we have
 $$A\lb zB(z)\rb+1=B(z)\qquad\text
{or}\qquad B\lb \frac z{1+A(z)}\rb =1+A(z)$$
and
 $$C\lb zB(z)\rb\cdot D(z)=(D(z)-1)\cdot B(z).$$
\endproclaim
\demo{Proof}
We only have to show the relation between $B(z)$, $C(z)$, and $D(z)$.
The crucial relation is the definition of the c-free cumulants
$R_k=R_k(\mu,\nu)$ by
 $$m_n(\mu)=\sum_{k=1}^n\sum\Sb l(1),\dots,l(k)\geq 0\\ l(1)+\dots+l(k)=n-k
\endSb R_k\cdot m_{l(1)}(\nu)\cdot\dots\cdot m_{l(k-1)}(\nu)\cdot
m_{l(k)}(\mu).$$
Now define
 $$\hat C(z):=\frac 1z C(z)=\sum_{n=1}^\infty R_nz^{n-1}.$$
Then we have
 $$\align
\hat C&\lb zB(z)\rb\cdot D(z)=\sum_{k=1}^\infty R_k\bigl(\sum_{l=0}^\infty
m_l(\nu)z^l\bigr)^{k-1}\bigl(\sum_{l(k)=0}^\infty m_{l(k)}(\mu)
z^{l(k)}\bigr) z^{k-1}\\
&=\sum_{k=1}^\infty R_k\sum_{l(1),\dots,l(k)\geq 0} m_{l(1)}(\nu)\cdot\dots
\cdot m_{l(k-1)}(\nu)\cdot m_{l(k)}(\mu)\cdot z^{l(1)+\dots+l(k)+(k-1)}\\
&=\sum_{n=1}^\infty\sum_{k=1}^n \sum\Sb l(1),\dots,l(k)\geq 0\\ l(1)+\dots
+l(k)=n-k\endSb R_k\cdot m_{l(1)}(\nu)\cdot\dots\cdot m_{l(k-1)}(\nu)\cdot
m_{l(k)}(\mu)\cdot z^{n-1}\\
&=\sum_{n=1}^\infty m_n(\mu)z^{n-1}\\
&=\frac 1z(D(z)-1),\endalign$$
hence
 $$\frac 1{zB(z)} C\lb zB(z)\rb\cdot D(z)=\frac 1z (D(z)-1),$$
which gives the assertion. \qed
\enddemo
Instead of dealing with the
generating power series $B(z)$ and $D(z)$ in the moments
it is usually more convenient to replace them by the corresponding
Cauchy-transforms
 $$g(z)=1/z\cdot B(1/z)\qquad\text{and}\qquad
G(z)=1/z\cdot D(1/z).$$
If we also replace the series $A(z)$ and $C(z)$ by the
$r/R$-transforms
 $$r(z)=A(z)/z\qquad\text{and}\qquad R(z)=\hat C(z)=C(z)/z,$$
then our main result can be rewritten as follows.
\proclaim{Theorem 5.2}
With the notations as above we have
 $$g(z)=\frac 1{z-r\lb g(z)\rb}
\qquad\text{or}\qquad g\lb r(z)+z^{-1}\rb=z$$
and
 $$G(z)=\frac 1{z-R\lb  g(z)\rb}.$$
\endproclaim
\demo{Examples}
1) Gaussian distribution as in Theorem 4.3.\newline
We have
 $$r(z)=\beta^2 z\qquad\text{and}\qquad R(z)=\alpha^2z,$$
which gives
 $$g(z)=\frac 1{z-\beta^2g(z)}
\qquad\text{and}\qquad
G(z)=\frac 1{z-\alpha^2g(z)},$$
which agrees with our calculations in Sect. 4.
Note that in our proof of Theorem 4.3 we used other combinatorial
identities than here. Our current machinery does not reproduce the proof
of 4.3, but it specializes to the formulas given in the remark after 4.3.
\newline
2) Poisson distribution as in Theorem 4.4.\newline
We have
 $$r(z)=\beta\frac 1{1-z}\qquad\text{and}\qquad R(z)=\alpha\frac 1{1-z},$$
which gives
 $$g(z)=\frac1{z-\beta\frac 1{1-g(z)}}
\qquad\text{or}\qquad
g(z)=\frac 1{z+(1-\beta)-zg(z)}$$
and
 $$G(z)=\frac 1{z-\alpha\frac 1{1-g(z)}}
\qquad\text{or}\qquad
G(z)=\frac 1{z+\frac \alpha\beta(1-\beta)-\frac\alpha\beta zg(z)},$$
in agreement with our calculations in Sect. 4.
\enddemo
\demo{Remark}
In the case of the boolean convolution
\hbox{$(\mu,\delta_0)=(\mu_1,\delta_0)\splus(\mu_2,\delta_0)$},
\linebreak which we mentioned in Remark 2 in Sect. 3, we have
$g(z)=g_{\delta_0}(z)=1/z$ and our formula in Theorem 5.2 reduces to
 $$G(z)=\frac 1{z-K(z)} \qquad\text{with}\qquad K(z)=R(1/z).$$
This simple formula reflects the simple structure of the underlying lattice
of interval partitions and offers the possibility for a far reaching
analytic treatment of the boolean convolution, in this respect see
\l Wor,SpW\r.
\enddemo
\vskip1cm
\heading
{\bf Acknowledgements}
\endheading
This work has been supported by
Polish National Grant, KBN 4019 (M.B.)
and a fellowship from the Deutsche
Forschungsgemeinschaft (R.S.).
\vskip1cm
\heading
{\bf References}
\endheading
\roster
\item"\l AGl\r"
Achieser, N.I., Glasmann, I.M.: {\it Theorie der linearen Operatoren im
Hilbert Raum}, Akademie-Verlag, Berlin 1960
\item"\l Bia\r"
Biane, P.: Some properties of crossings and partitions. Preprint, Paris 1994
\item"\l Boz1\r"
Bo\D zejko, M.: Positive definite functions on the free group and the
non commutative Riesz product. Boll. Un. Mat. Ital. 5A,
13-21 (1986)
\item"\l Boz2\r"
Bo\D zejko, M.: Uniformly bounded representations of free groups.
J. Reine Angew. Math. 377, 170-186 (1987)
\item"\l BSp\r"
Bo\D zejko, M., Speicher, R.: $\psi$-independent and symmetrized white noises.
In {\it Quantum Probability and Related Topics VI}
(ed. L. Accardi), World Scientific, Singapore 1991, 219-236.
\item"\l CDu\r"
Cabanal-Duvillard, T.: Variation quantique sur l'independance:
La $\sigma$-independance. Preprint, Paris 1993
\item"\l Ede1\r"
Edelman, P.H.: Chain Enumeration and Non-Crossing Partitions. Discr. Math. 31,
171-180 (1980)
\item"\l Ede2\r"
Edelman, P.H.: Multichains, Non-Crossing Partitions and Trees. Discr. Math.
40, 171-179 (1982)
\item"\l Gir\r"
Girko, W.L.: {\it Random Matrices}, Kiev 1975 (in Russian)
\item"\l Haa\r"
Haagerup, U.: An Example of a Non Nuclear $C^*$-Algebra which has the
Metric Approximation Property. Invent. Math. 50, 279-293 (1979)
\item"\l HiP\r"
Hilton, P., Pederson, J.: Catalan Numbers, Their Generalization, and Their
Uses. Math. Intelligencer 13, No. 2, 64-75 (1991)
\item"\l Kre\r"
Kreweras, G.: Sur les partitions non croisees d'un cycle. Discr. Math. 1,
333-350 (1972)
\item"\l Maa\r"
Maassen, H.: Addition of freely independent random variables. J. Funct.
Anal. 106, 409-438 (1992)
\item"\l Nic\r"
Nica, A.: Crossings and embracings of set-partitions and $q$-analogues of
the logarithm of the Fourier transform. Preprint, Berkeley 1994
\item"\l NSp1\r"
Neu, P., Speicher, R.: A self-consistent master equation and a new
kind of cumulants. Z. Phys. B 92, 399-407 (1993)
\item"\l NSp2\r"
Neu, P., Speicher, R.: Non-linear master equation and non-crossing cumulants.
Preprint HD-AM-NSP-93-1. To appear in {\it Quantum Probability and
Related Topics IX} (ed. L. Accardi), World Scientific, Singapore
\item"\l Pou\r"
Poupard, Y.: Etude et denombrement paralleles des partitions non croisees
d'un cycle et coupage d'un polygone convexe. Discr. Math. 2, 279-288 (1972)
\item"\l Sim\r"
Simion, R.: Combinatorial statistics on non-crossing partitions. Preprint,
Washington 1991
\item"\l SiU\r"
Simion, R., Ullman, D.: On the structure of the lattice of
non-crossing partitions. Discr. Math. 98, 193-206 (1991)
\item"\l Spe1\r"
Speicher, R.: A New Example of \lq Independence' and \lq White Noise'.
Probab. Th. Rel. Fields 84, 141-159 (1990)
\item"\l Spe2\r"
Speicher, R.: Multiplicative functions on the lattice of non-crossing
partitions and free convolution.
Math. Ann. 298, 611-628 (1994)
\item"\l SpW\r"
Speicher, R., Woroudi, R.: Boolean convolution. Preprint, Heidelberg 1993
\item"\l Voi1\r"
Voiculescu, D.: Symmetries of some reduced  free product
$C^*$-algebras. In
{\it Operator Algebras and their Connection with
Topology and Ergodic Theory},
LNM 1132, Springer, Heidelberg 1985, 556-588
\item"\l Voi2\r"
Voiculescu, D.: Addition of certain non-commuting random variables.
J. Funct. Anal. 66, 323-346 (1986)
\item"\l VDN\r"
Voiculescu, D., Dykema, K., Nica, A.: {\it Free Random Variables},
AMS 1992
\item"\l vWa\r"
von Waldenfels, W.: An approach to the theory of pressure
broadening of spectral lines. In: Lecture Notes in Math. 296,
Springer, Heidelberg 1973, 19-69
\item"\l Wor\r"
Woroudi, R.: Boolesche Faltung. Diplomarbeit, Heidelberg 1993
\endroster
\vskip1cm
\heading
{\bf Figure Captions}
\endheading
{\bf Fig. 1:}
c-free Gaussian distribution $\nu_{\alpha,\beta}$
for fixed $\beta=1$ and six different values of $\alpha$;
vertical double lines indicate $\delta$-peaks
\newline
{\bf Fig. 2:}
c-free Poisson distribution $\pi_{\alpha,\beta}$ for fixed $\alpha=1$
and six different values of $\beta$; vertical double lines indicate
$\delta$-peaks; note that the $\delta$-peak at $z_0$ lies first
on the right side of the continuous spectrum, then it dissapears and
reappears again on the left side of the continuous spectrum
\enddocument